\documentclass[twoside]{irmaems}
\usepackage{amssymb} %use \usepackage[psamsfonts]{amssymb} if your fonts are from Y&Y/Blue Sky Research
\usepackage{amsmath} %use \usepackage{amsmath,cmmib57} if your fonts are from Y&Y/Blue Sky Research
\usepackage{latexsym}

\usepackage{graphicx}

\renewcommand\theequation{\arabic{section}.\arabic{equation}}

\theoremstyle{definition} %%% for statements in roman typeface

 \newtheorem{definition}{Definition}[section]
 \newtheorem{remark}[definition]{Remark}

\theoremstyle{plain}      %%% for statements in italic typeface

 \newtheorem{theorem}[definition]{Theorem}

\title{Hilbert's Fourth Problem}

\begin{document}

\title{On Hilbert's fourth problem}

\author{Athanase Papadopoulos}\address{Institut de Recherche Math\'ematique Avanc\'ee,\\
Universit{\'e} de Strasbourg and CNRS,\\
7 rue Ren\'e Descartes,\\
 67084 Strasbourg Cedex, France.\\
 email : papadop@math.unistra.fr} 

\date{\today}

\maketitle

 \begin{abstract} Hilbert's fourth problem asks for the construction and the study of metrics on subsets of projective space for which the projective line segments are geodesics. Several solutions of the problem were given so far, depending on more precise interpretations of this problem, with various additional conditions satisfied. The most interesting solutions are probably those inspired from an integral formula that was first introduced in this theory by Herbert Busemann. Besides that, Busemann and his school made a thorough investigation of metrics defined on subsets of projective space for which the projective lines are geodesics and they obtained several results, characterizing several classes of such metrics.  We review some of the developments and important results related to Hilbert's problem, especially those that arose from Busemann's work, mentioning recent results and connections with several branches of mathematics, including Riemannian geometry, the foundations of mathematics, the calculus of variations, metric geometry and Finsler geometry.  
  \end{abstract}


\noindent 
 \emph{Acknowledgements : }
 I would like to thank Juan Carlos \'Alvarez Paiva, Charalampos Charitos, Jeremy Gray, Fran\c cois Laudenbach,  Marc Troyanov and Sumio Yamada for comments on a preliminary version of this paper.
 This work was partially supported by  the French ANR project FINSLER.
 
  \medskip
 \noindent  The final version of this paper will appear in the \emph{Handbook of Hilbert geometry}, published by the European Mathematical Society (ed. A. Papadopoulos and M. Troyanov).

  \vfill\eject

\section{Introduction}

On the 8th of August 1900, at the Second International Congress of Mathematicians held in Paris, David Hilbert delivered a lecture titled ``The future problems of mathematics", in which he presented a collection of open problems. Hilbert was forty-eight years old and was considered as one of the leading and most versatile mathematicians of his time. His list of problems turned out to be a working document for a large part of the mathematical research that was conducted in  the twentieth century and it had a great influence on several generations of mathematicians. A first short paper based on that lecture appeared in 1900 in the newly founded Swiss journal \emph{L'Enseignement Math\'ematique} \cite{Hilbert1900}, and after that several written versions of unequal lengths were published (see \cite{Grattan} for a commentary on these versions). This includes a  paper published in the \emph{Bulletin of the American Mathematical Society} (\cite{Hilbert-Problems}, 1901), containing a commented set of twenty-three problems.\footnote{In the mid 1990s, a twenty-fourth  problem was discovered in Hilbert's massive files in G\"ottingen. This  problem asks for ``\emph{simpler} proofs and criteria for simplicity". The interested reader can refer to the reports in \cite{Thiele-Wos} and  \cite{Thiele} on that intriguing problem that had remained hidden for almost a century. We also note that at the Paris 1900 congress, Hilbert presented only 10 problems. C. Reid, who recounts the story in \cite{Ried}, says (on p. 81): ``In an effort to shorten his talk as Minkowski and Hurwitz had urged, Hilbert presented only ten problems out of a total of 23 which he had listed in his manuscript".}  Today, several of these problems are solved, but there is still very active research around them (the solved ones as well as the unsolved). 

 As a general  introduction to Hilbert problems, we refer the reader to the two semi-popular books that appeared on the occasion of the anniversary of these problems:  J. Gray's  \emph{The Hilbert Challenge} \cite{Gray2000} and B. H. Yandell's \emph{The Honors Class: Hilbert's Problems and Their Solvers} \cite{Yandell2002}. The books \cite{Browder} and \cite{Alexandrov} contain collections of more technical articles on the developments of the Hilbert problems.

 Hilbert geometry, which is the topic of this Handbook, is closely connected with Hilbert's Problem IV. The problem is titled  ``The problem of the straight line as the shortest distance between two points". There are several interpretations of this problem and we shall discuss them below. Soon after Hilbert formulated it, the problem was reduced by G. Hamel -- who was one of Hilbert's student -- to a question on metrics on \emph{convex} subsets of Euclidean (or affine) spaces. If $\Omega$ is a convex subset of a Euclidean space, then the restriction to $\Omega$ of the ambient Euclidean metric is a metric on $\Omega$ satisfying Hilbert's requirement, that is, the (restriction to $\Omega$ of the) Euclidean straight lines are geodesics for the restriction to $\Omega$ of the Euclidean metric. One form of Hilbert's problem asks for a characterization of \emph{all} metrics on $\Omega$ for which the Euclidean lines in $\Omega$  are geodesics. A natural question is then to find a characterization of metrics on $\Omega$ for which the Euclidean lines are the unique geodesics. One may put additional requirements on the metrics and try to characterize some special classes of metrics satisfying Hilbert's requirements, such as \emph{geodesically complete} metrics, that is, metrics whose geodesics (parametrized by arclength) can be extended infinitely from both sides.  The Hilbert metrics are examples of geodesically complete metrics satisfying Hilbert's problem requirement.
 
 In what follows, I shall present some ideas and make comments, several of them having a historical character, about Hilbert's Problem IV and its ramifications, and I will concentrate on the ideas that arose from the work of Busemann. In particular, I generally followed the rule of not stating in the form \emph{Theorem: etc.} anything that is not found in Busemann's work.

\section{Precursors}

Hilbert's problems did not come out of nowhere, and several of them have a history.
Problem IV, like most of the others, arose as a natural question at that time. In the French version of the statement (see the appendix at the end of the present paper), Hilbert mentions related work of Darboux. Darboux, in the work that Hilbert refers to, mentions Beltrami. We shall review some of Beltrami's work related to the subject, but before that, we report on some works of Cayley and Klein.

In a 1859 paper called \emph{A sixth Memoir upon Quantics} \cite{Cayley1859}, Cayley
\footnote{Arthur Cayley (1821-1895) grew up in a family of English merchants settled in Saint-Petersburg, and he lived in Russia until the age eight, when his family returned to England. Cayley is one of the main inventors of the theory of invariants. These include invariants of algebraic forms (the determinant being an example), but also algebraic invariants of geometric structures and the relations they satisfy (``syzygies"). Cayley studied both law and mathematics at Trinity College, Cambridge, from where he graduated in 1842. was talented in both subjects and as an undergraduate he wrote several papers, three of which were published in the \emph{Cambridge Mathematical Journal}. The subject of these papers included determinants, which became later on one of his favorite subjects of study. After a four-year position at Cambridge university, during which he wrote 28 papers for the Cambridge journal, did not find any job to continue in academics. Looking for a new profession, law was naturally his second choice. He worked as a lawyer for 14 years, but he remained interested in mathematicians, and he wrote about 250 mathematical papers. Cayley was appointed professor of mathematics at Cambridge in 1863. His list of papers includes about 900 items, on all the fields of mathematics of his epoch. The first definition of abstract group is attributed to him, cf. his paper \cite{Cayley-group}. Cayley is also one of the first discoverers of geometry in $n$ dimensions. In his review of Cayley's \emph{Collected Mathematical Papers} edition in 13 volumes, G. B. Halsted writes: ```Cayley not only made additions to every important subject of pure mathematics, but whole new subjects, now of the most importance, owe their existence to him. It is said that he is actually now the author most frequently quoted in the living world of mathematicians" \cite{Halsted-Cayley}. The reader is referred to the interesting biography by Crilly \cite{Crilly}.} gave an interesting definition of a \emph{projective metric on the plane}, that is, a metric on the projective plane such that the projective lines are geodesics for that metric. More specifically, Cayley started with a conic (which he called the \emph{absolute}) in the projective plane, and he associated to it a metric $d$ such that for any three points $P_1, P_2, P_3$ in that order on a projective line, we have $d(P_1,P_2)+ d(P_2,P_3)= d(P_1,P_3)$. A few years later, Klein, in his two papers \cite{Klein-Ueber}  and \cite{Klein-Ueber1}, re-interpreted Cayley's ideas using the notion of cross-ratio, and he showed that by choosing the absolute to be a real conic in the projective plane and using Cayley's construction, the interior of that conic is a model of hyperbolic geometry. Beltrami, in his paper \cite{Beltrami1868} had already noticed this model, namely, he showed that the unit disc in the plane, with the Euclidean chords taken as geodesics, is a model of the Lobachevsky plane with its non-Euclidean geodesics, but he did not have the formula for the distance in terms of the cross ratio. Klein's model, with the cross ratio formula, constitutes what is called today the \emph{Klein model}, or the  \emph{Beltrami-Klein model}, or also the  \emph{Beltralmi-Cayley-Klein model} of hyperbolic space.

It is worth noting that Klein found also models of the elliptic and of the Euclidean plane within the projective plane that use the cross ratio and a conic at infinity. In the case of the Euclidean plane, the conic is imaginary and the construction involves complex numbers. We refer the reader to Klein's paper \cite{Klein-Ueber}, and also to the commentary \cite{ACP1}.

It is may also be worth recalling here that Cayley considered the use of the cross ratio by Klein as a sort of circular reasoning, because he thought that the definition of the cross ratio is dependent upon some underlying Euclidean geometry, whereas Klein, in his paper \cite{Klein-Ueber}, was aiming for a  definition of the hyperbolic metric (and also, Euclidean and spherical metrics) based only on projective geometry.  We can quote here Cayley from his comments on his paper \cite{Cayley1859} in Volume II of  his
   \emph{Collected mathematical papers} edition \cite{Cayley-collected} (p. 605):
 \begin{quote}\small
  I may refer also to the memoir, Sir R. S. Ball ``On the theory of content," \emph{Trans. R. Irish Acad.} vol. {\sc xxix} (1889), pp. 123--182, where the same difficulty is discussed. The opening sentences are -- ``In that theory [Non-Euclidian geometry] it seems as if we try to replace our ordinary notion of distance between two points by the logarithm of a certain anharmonic ratio. But this ratio itself involves the notion of distance measured in the ordinary way. How then can we supersede the old notion of distance by the non-Euclidian notion, inasmuch as the very definition of the latter involves the former ?"
\end{quote}
In the same line of thought, we quote Genocchi, an influential and well-established Italian mathematician of the second half of the nineteenth century\footnote{Angelo Genocchi (1817-1889) was in Italian mathematician who was also a very active politician. Like Cayley, he worked for several years as a lawyer, and he taught law at the University of Piacenza, the city where he was born, and at the same time he was cultivating mathematics with passion. In 1859, Genocchi was appointed professor of mathematics at the University of Torino, and he remained there until 1886. During the academic year 1881-82, Guiseppe Peano served as his assistant, and he helped him later on  with his teaching, after Genocchi became disabled after of an accident. During several years, Genocchi was the main specialist in number theory in Italy.}  who wrote (\cite{Genocchi} p. 385): 
\begin{quote}\small
From the geometric point of view, the spirit may be chocked by certain definitions adopted by Mr. Klein: the notions of distance and angle, which are so simple, are replaced by complicated definitions [...] The statements are extravagant." 
\end{quote}
This is just to show that some of the most prominent mathematicians were missing an important idea, namely that the cross ratio does not depend on the underlying Euclidean geometry.

The second paper of Klein \cite{Klein-Ueber1} contains an elaboration on some of the ideas contained in the first one. In the meanwhile Beltrami found two models of hyperbolic geometry, and Klein noticed that his model coincides with one of the models that Beltrami found. This is the origin of the so-called Klein-Beltrami model of hyperbolic space. Klein's aim in the two papers that we mentioned was to include the three constant-curvature geometries in the setting of projective geometry (Cayley had done this only for Euclidean geometry), and although Klein did not formulate it in these words, he obtained models for these three geometries in which lines in projective space are geodesics. This fact was highlighted by Beltrami, who obtained in \cite{Beltrami1868} the disc-model of hyperbolic geometry, where the Euclidean straight lines are the geodesics (but Beltrami did not have Klein's formula involving the cross-ratio.

 One may add here that two years before Klein's work on the subject, Beltrami wrote the following to Hou\"el\footnote{Guillaume Jules Ho\"uel (1823--1886) is a major figure in the history of non-Euclidean geometry. There is a very interesting correspondence between him and Beltrami, and sixty-five of these letters were edited by Boi, Giacardi and Tazzioli in 1998  \cite{Beltrami-Boi}. We refer to  that correspondence several times in the present paper. Ho\"uel  translated  into French and published in French and Italian journals works by several other authors on non-Euclidean geometry, including Bolyai, Beltrami, Helmholtz, Riemann and Battaglini.  In his first letter to Ho\"uel, written on November 18, 1868, Beltrami (who was seven years younger than Ho\"uel), starts as follows: ``Professor, I take the liberty of sending you, in a file, two memoirs, the object of the first one is a real   construction of non-euclidean geometry, and the other, much less recent, contains the proof of some analytic results upon which this construction relies, but which, apart from that, have no immediate relation with that question". The two memoirs are the  \emph{Saggio di Interpretazione della geometria non-Euclidea} \cite{Beltrami1868} and the \emph{Risoluzione del problema: Riportare i punti di una superficie sopra un piano in modo che le linee geodetiche vengano rappresentate da linee rette} \cite{Beltrami1865}. It seems that Ho\"uel was the first to notice that Beltrami's work on the pseudo-sphere implies that the parallel postulate is not a consequence of the other Euclidean postulates. Beltrami, in a letter dated January 2, 1870 (\cite{Beltrami-Boi} p. 114), writes: ``Your conclusion on the impossibility of proving in the plane the eleventh axiom, as far as it follows from my work, is also mine". See the discussion in \cite{Henry-Nabonnand-H}, \S 1.3 of the Introduction.} (letter dated July 29, 1869 \cite{Beltrami-Boi} p. 96-97):
   \begin{quote}\small
   The second thing [I plan to do in addition to the memoir I sent you]\footnote{Beltrami, in this letter, gave a favorable response to Ho\"uel who had proposed to translate his two memoirs \emph{Saggio di Interpretazione della geometria non-Euclidea} \cite{Beltrami1868} and \emph{Teoria fondamentale degli spazi di curvatura costante} \cite{Beltrami-Teoria} into French and to publish the translations in the \emph{Annales de l'\'Ecole Normale Sup\'erieure}, and he is talking to him here about two new projects he has.} will be the most important, if I succeed to give it a concrete form, because up to now it only exists in my head in the state of a vague conception, although without any doubt it is based on the truth. This is the conjecture of a close analogy, and may be of an identity, between pseudo-spherical geometry\footnote{This is the term used by Beltrami to denote hyperbolic geometry.} and the theory of Mr Cayley on the \emph{analytical origin of metric ratios}, using  an \emph{absolute} conic (or quadric). I almost did not know anything about that theory, when I was taken by the identity of certain forms. However, since the theory of invariants plays there a rather significant role and because I lost sight of it since a few years now, I want to do it again after some preliminary studies, before I address this comparison.
      \end{quote}
   In another letter to Ho\"uel, written on July 5, 1872, Beltrami regrets the fact that he let Klein outstrip him (\cite{Beltrami-Boi} p. 165):
      \begin{quote}\small
The principle which has directed my analysis\footnote{Beltrami refers here to a note \cite{beltrami-Osservatione} which he had just published in the \emph{Annali di Matematica}.} is precisely that which Mr Klein has just developed in his recent memoir\footnote{Beltrami refers to Klein's paper \cite{Klein-Ueber}.} on non-Euclidean  geometry, for 2-dimensional spaces. In other words, from the analytic point of view, the geometry of spaces of constant curvature is nothing else than Cayley's doctrine of the absolute. I regret very much to have let Mr Klein supersede me on that point, on which I had already assembled some material, and it was my mistake for not giving to this matter enough importance.            \end{quote}
There are comments on this work of Klein in the paper \cite{ACP1}.

 We next talk about Beltrami's work. 
 
In 1865, Beltrami wrote a paper \cite{Beltrami1865} whose title is \emph{Risoluzione del problema: Riportare i punti di una superficie sopra un piano in modo che le linee geodetiche vengano rappresentate da linee rette} (A solution of the problem: Transfer the points of a surface over a plane in such a way that the geodesic lines are represented by straight lines). From the title, it is clear that the problem bears a direct relation to Hilbert's problem.  Indeed, if such a transfer map exists, then, by pushing forward the metric of the given surface to a metric on the Euclidean plane, we get a metric on that plane for which the Euclidean straight lines are geodesics.

Beltrami's work on this subject was in the framework of the differential geometry of surfaces in the style developed by Gauss. In the introduction to his paper \cite{Beltrami1865}, Beltrami says that the major part of the research done before him on similar questions was concerned with the question of conservation of angles or of area proposition, and that even though these two properties are regarded as the simplest and most important properties of such transfer maps, there are other properties that one might want to preserve. Indeed, the kind of question in which Beltrami was interested was motivated by the theory of drawing of geographical maps. Beltrami declares that since the projection maps that are used in this science were principally concerned in the measure of distances, one would like to exclude projection maps where the images of distance-minimizing curves are too remote from straight lines. He mentions in passing that the central projection of the sphere is the only map that transforms geodesics into straight lines. As a matter of fact, this is the question that motivated Beltrami's research on the problem referred to in the title of his essay,  ``Transfer the points of a surface over a plane etc.". Beltrami then writes that beyond its applications to geographic map drawing, the solution of the problem may lead to ``a new method of geodesic calculus, in which the questions concerning geodesic triangles on surfaces can all be reduced to simple questions of plane trigonometry". He finally acknowledges that his investigations on this question did not lead him to any general solution, and that the case of the central projection of the sphere is the only one he found where the required condition is realized.

Beltrami proved an important local result in that paper, viz. he showed that the only surfaces that may be (locally) mapped to the Euclidean plane in such a way that the geodesics are sent to Euclidean lines must have constant Gaussian curvature.

In fact,  Beltrami starts with a surface whose line element is given, in the tradition of Gauss,  by 
\[ds^2= Edu^2+2Fdudv +Gdv^2,\]
and after several pages of calculations he ends up with the following values for the functions $E$, $F$ and $G$:
\begin{equation}\label{Beltram-f}
\left\{
     \begin{array}{lll}
        E   = \frac{R^2(v^2+a^2)}{(u^2+v^2+a^2)^2}  \\
       F = \frac{-R^2uv}{(u^2+v^2+a^2)^2} \\
 G= \frac{R^2(u^2+a^2)}{(u^2+v^2+a^2)^2}
     \end{array}
  \right.
\end{equation}
Beltrami then uses the following formula for the curvature, which he derived in his \emph{Ricerche} \cite{Beltrami-Ricerche}, art. XXIV:
  \[
  \frac{1}{R_1R_2}= -\frac{1}{2\sqrt{EG-F^2}}\left( \frac{\partial}{\partial u}\left(
  \frac{\frac{\partial G}{\partial u}-\frac{F}{E}\frac{\partial E}{\partial v}}{\sqrt{EG-F^2}}\right)
  +\frac{\partial}{\partial v}\left(
  \frac{\frac{\partial E}{\partial v}-2\frac{\partial F}{\partial u}  +\frac{F}{E}\frac{\partial E}{\partial u}}{\sqrt{EG-F^2}} \right) \right)
  \]
 After a computation, Beltrami finds that the curvature is constant and equal to
  \[  
  \frac{1}{R_1R_2}=\frac{1}{R^2}.
  \]
  He concludes from this:  ``Therefore our surfaces are those of constant curvature. In particular, if the quantity $R$ is real, the formulae (\ref{Beltram-f}) are responsible for all surfaces that are applicable on the sphere of radius $R$".
We recall that it is a consequence of Gauss's \emph{Theorema egregium} that a surface of constant positive curvature is locally isometric to a sphere.    
 
 Beltrami summarizes his result in the following statement:
 
\begin{theorem}
 The only surfaces that can be represented over a plane, in such a way that to every point corresponds a point and to every geodesic line a straight line are those whose curvature is everywhere constant (positive, negative or zero). When this constant curvature is zero, the correspondence does not differ from the ordinary homography. When it is nonzero, this correspondence is reducible to the central projection over the sphere and to its homographic transformations.  
 \end{theorem}
Beltrami adds:
 \begin{quote}\small
 Since among all the surfaces of constant curvature the only one that can have applications in the theory of geographic maps and in geodesy is probably the spherical surface, in this way from the point of view of these applications, what we asserted confirms that the only solution to the problem is obtained essentially by the central projection. 
  \end{quote}

 Thus, in his conclusion, Beltrami considered only the case where the curvature is positive, and therefore the only surface he found that has the required property (apart from the plane itself) is the sphere (of a certain radius).  We also refer the reader to the chapter \cite{Tro-Finsler} in this Handbook for an exposition with a detailed proof of this result of Beltrami.

 Two years later,  Beltrami wrote another paper, \emph{Saggio di Interpretazione della geometria non-Euclidea} \cite{Beltrami1868}, which soon became famous, in which he worked on surfaces of constant negative curvature, a context in which he gave an interpretation of Lobachevsky's geometry.  In that  paper, Beltrami defined a family of metrics on the disc of radius $a>0$ centered at the origin of the Euclidean plane, using the formula
\begin{equation} \label{BN1}
ds^2=R^2\frac{(a^2-v^2)du^2+2ududv + (a^2-u^2)dv^2}{(a^2-u^2-v^2)^2},
\end{equation}
where $ds^2$ denotes as usual the square of the infinitesimal length element. He showed that the Gaussian curvature of that surface is constant and equal to $1/R^2$, and that its geodesics are the Euclidean straight lines. He also stated the following:
\begin{quote}\small
It follows from what precedes that the geodesics of the surface are represented in their total (real) development by the chords of the limit circle, whereas the prolongation of these chords in the exterior of the same circle are devoid of any (real) representation.\footnote{The English translation  in \cite{Sources} of this sentence from \cite{Beltrami1868} is not very faithful.}
\end{quote}

 Beltrami gave in the paper \cite{Beltrami1868} the first Euclidean model of hyperbolic plane (the so-called the ``Beltrami" model, the ``Klein" model, the ``Klein-Beltrami" model or the ``Cayley-Klein-Beltrami" model\footnote{The four names are reasonable, see e.g. the comments in \cite{ACP1}}), and at the same time he gives a new example of a metric on the plane satisfying Hilbert's problem. We mention by the way that this research by Beltrami was motivated by his reading of the works of Lobachevsky. In a letter to Ho\"uel, written on November 18, 1868, Beltrami writes: ``That writing [The \emph{Saggio} etc.] was written during last year's fall, after thoughts that abounded in me at the epoch of the publication of your translation of Lobatschewsky. I think that the idea of constructing non-Euclidean geometry on a perfectly real surface was entirely new"
 (\cite{BGT}) p. 65-66).

In Note I at the end of the paper \cite{Beltrami1868}, Beltrami recalls the metric of constant positive curvature $1/R^2$ that he defined in his 1865 paper \cite{Beltrami1865}, whose square of the infinitesimal line element is given by the formula
\begin{equation} \label{BN2}
ds^2=R^2\frac{(a^2+v^2)du^2-2ududv + (a^2+u^2)dv^2}{(a^2+u^2+v^2)^2}
\end{equation}
and he makes an interesting remark on the relation between the two metrics (\ref{BN1}) and (\ref{BN2}), namely, that the formula for one metric can be obtained from the other one by replacing the constants $R$ and $a$ by $R\sqrt{-1}$ and $a\sqrt{-1}$, thus confirming the statement that hyperbolic geometry is in some sense spherical geometry worked out on a ``sphere of imaginary radius".

A modern proof of some of Beltrami's results are contained in Busemann's book \cite{Busemann1955}, Chapter II, \S 15. Busemann's introduces this result as follows:
\begin{quote}\small
Although the methods are quite foreign to the rest of the book, we prove this fundamental \emph{Theorem of Beltrami} in Section 15, because it is by far the most striking example of a Riemannian theorem without a simple analogue in more general spaces.
\end{quote}
Busemann states Beltrami's result as follows (see \cite{Busemann1955} p. 85, where Busemann refers to Blaschke \cite{Blaschke1936}): 
\begin{theorem}
Let $C$ be a connected open set of the projective plane $\mathbb{P}^2$ equipped with a Riemannian metric whose geodesics lie on projective lines. Then the Gauss curvature is constant. 
\end{theorem}
 
 For dimension $\geq 3$, Busemann refers to Cartan \cite{Cartan1928}.
 
 \begin{remark}[I owe this remark to Juan Carlos \'Alvarez Paiva]
It is not sufficiently known that Beltrami's theorem holds with very weak assumptions on the regularity on the metric. This was first noticed by Hartman and Wintner in an interesting paper \cite{HW}, in which they prove that if a Riemannian metric is  continuous and projective, then it is real analytic. Pogorelov rediscovered the result of Hartman and Wintner several years later.

 \end{remark}

Darboux, who is quoted by Hilbert in his formulation of his fourth problem, in his \emph{Le\c cons sur la th\'eorie g\'en\'erale des surfaces} (Troisi\`eme partie), following Beltrami's work, considers, in the setting of the calculus of variations, a similar problem. He writes: ``Proposons-nous d'abord le probl\`eme de M. Beltrami: \emph{Trouver toutes les surfaces qui peuvent \^etre repr\'esent\'ees g\'eod\'esiquement sur le plan}." (\cite{Darb} p. 59). That is to say, Darboux considers Beltrami's problem of sending a surface into the plane by a map that sends the surface geodesics to the plane straight lines. According to  Busemann (\cite{Busemann1970} p. 56),  Darboux proposed several questions that are related to Hilbert's problem, some of them being more general than the one in hand.\footnote{I owe the following not to Juan Carlos \'Alvarez Paiva: Darboux's ``solution" of Hilbert's fourth problem (before it was even stated) misses however one of the most 
important features of the problem: by considering the variational problem for integrands of the form $f(x,y,y')$ instead of integrands of the form $L(x,y, x', y')$ that are homogeneous of degree one in the velocities $x'$ and $y'$, Darboux gives many variational problems whose extremals are straight lines, but which cannot define metrics because they are not even defined on the whole tangent bundle (minus the origin). Darboux's approach can be modified and be made to yield yet another solution to Hilbert's fourth problem in dimension two in the smooth asymmetric case. Gautier and \'Alvarez-Paiva did this, and later found out that  it is also implicit in Hamel's paper \cite{Hamel1903}.} In fact, one can relate the ``problem of geographical maps", that we already mentioned to problems that were studied by Gauss, Lambert, Euler, Lagrange, Darboux, Liouville and Bonnet. Gauss, for instance, formulated the problem of ``representing a surface on the plane so that the representation is similar, at the very small level, to the original surface". This means that we ask for the similarity at the infinitesimal level for lines and angles between them, see Gauss's quotes and the comments in the paper \cite{Nabonnand}. Using a more modern language, this is the question of finding conformal coordinates.

 \section{Hilbert's problem IV}\label{s:H4}
  
 The best way in which one can have an idea on what Hilbert's Problem IV is to read Hilbert's text, and this is why we reproduced it in the appendix at the end of this paper.   It is also good to remember and to meditate on Hilbert's conclusion to his \emph{Mathe\-matical Problems}.\footnote{The translation is from the article in the 1902 Bulletin of the AMS. \cite{Hilbert-Problems}}:
 
\begin{quote}\small
The problems mentioned are merely samples of problems, yet they will suffice to show how rich, how manifold and how extensive the mathematical science of today is, and the question is urged upon us whether mathematics is doomed to the fate of those other sciences that have split up into separate branches, whose representatives scarcely understand one another and whose connection become ever more loose. I do not believe this nor wish it. Mathematical science is in my opinion an indivisible whole, an organism whose vitality is conditioned upon the connection of its parts. For with all the variety of mathematical knowledge, we are still clearly conscious of the similarity of the logical devices, the \emph{relationship} of the \emph{ideas} in mathematics as a whole and the numerous analogies in its different departments. We also notice that, the farther a mathematical theory is developed, the more harmoniously and uniformly does its construction proceed, and unsuspected relations are disclosed between hitherto separate branches of the science. So it happens that, with the extension of mathematics, its organic character is not lost but only manifests itself more clearly. 

But we ask, with the extension of mathematical knowledge, will it not finally become impossible for the single investigator to embrace all departments of this knowledge ? In answer let me point out how thoroughly it is ingrained in mathematical science that every real advance goes hand in hand with the invention of sharper tools and simpler methods which at the same time assist in understanding earlier theories and cast aside older more complicated developments. It is therefore possible for the individual investigator, when he makes these sharper tools and simpler methods his own, to find his way more easily in the various branches of mathematics than is possible in any other science.

The organic unity of mathematics is inherent in the nature of this science, for mathematics is the foundation of the exact knowledge of natural phenomena. That it may completely fulfill his high mission, may the new century bring it gifted masters and many zealous and enthusiastic disciples.
\end{quote}

 We now make some comments on Hilbert's Problem IV.

 The formulation by Hilbert and the scope of his fourth problem are very broad compared to other problems in his list,\footnote{\label{Problems}As examples of Hilbert problems that are more precisely stated, we mention Hilbert's Problem III, asking whether two 3-dimensional Euclidean polyhedra of the same volume are scissors-equivalent. This problem was solved in the negative by M. Dehn soon after Hilbert formulated it. We shall elaborate on this in Footnote \ref{n:Dehn} below. Another well-known example of a precise problem in Hilbert's list is Part 2 of Problem VII, which asks whether given any algebraic number  $a\not= 0,1$ and given any irrational algebraic number $b$,  $a^b$ is transcendental. This question was solved by A. O. Gelfond in 1934 and the proof was refined by T. Schneider in 1935. The result is called now the Gelfond-Schneider Theorem. 
 But there are also problems in Hilbert's list which are vaguely stated (and some of them are so because of their nature). For instance, Problem VI asks for a ``mathematical treatment of the axioms of physics".} and this problem admits several interpretations and generalizations. In fact, the problem was considered at several times in history as being solved.\footnote{One can search for papers titled ``A solution to Hilbert's Problem IV" and similar titles; there are several.} We shall recall in later sections how several mathematicians who worked on this problem interpreted it in their own manner and contributed to a solution. We mention right away the names of G. Hamel and P. Funk who worked on it in the early years of the 20th century, the fundamental work of H. Busemann over a span of 50 years (1930-1980), and the works of A. V. Pogorelov, Z. I. Szab\'o, R. Ambartzumian, R. Alexander, I. M. Gelfand, M. Smirnov and J.-C. \'Alvarez Paiva. We shall elaborate on some of these works below.
 
 There are at least five natural settings in which one can approach Hilbert's Problem IV: 
\begin{enumerate}
\item The setting of Riemannian geometry. In this setting, the problem was settled by Beltrami.
\item  The setting of the calculus of variations,  which was the one of Darboux in the work we mentioned. In this setting, Hilbert's problem is often described as an ``inverse problem"; that is, one looks for a Lagrangian (defining the metric) which is associated to a given set of extremals (which are the geodesic lines). As is well known, the methods of the calculus of variations started in the works of Johann and Jakob Bernoulli, of Euler and of Lagrange.
\item The setting of metric geometry, where one does not assume any differentiability. This setting is at basis of the methods of Busemann and of the large amount of interesting works that arose from them.
\item   The setting of Finsler geometry, where the metric is associated to a norm on each tangent space of the manifold. A minimum of differentiability is required in this setting, at least to give a meaning to tangent vectors. This setting can be considered as lying at the border of metric and of differential geometry.
\item  The setting of the foundations of geometry. This was probably the most important setting for Hilbert, and we shall elaborate on it below. In fact, this setting is also inherent in one way or another in all the different settings mentioned above.
\end{enumerate}
There are also other points of view and relations with other questions. For instance, R. Alexander established relations between  Hilbert's Problem IV and zonoid theory (see \cite{Alexander} and the paper \cite{Szabo} by Z. Szab\'o), S. Tabachnikov made relations with the theory of magnetic flows (see \cite{Tabachnikov}).  I. M. Gelfand, M. Smirnov and J.-C. \'Alvarez Paiva developed relations with symplectic geometry (see \cite{AGS}, \cite{Alvarez2000} and  \cite{Alvarez2005}). The result in 
\cite{Alvarez2000} can be stated as follows: There exists a (twistorial)
correspondence between smooth, symmetric solutions of Hilbert's fourth problem in three dimensional projective space and anti-self-dual symplectic forms in $S^2 \times S^2$. 

In any case, Hilbert's Problem IV and the various approaches to solve it had an important impact on the development of several fields of mathematics, including all the ones we mentioned, to which we can add convex geometry and integral geometry. The developments that arose from that problem brought together major ideas originating from the various fields mentioned.

Let us return to Hilbert's statement, ``the problem of the straight line as the shortest distance between two points". Today, when we talk about straight lines being the shortest distance, we think of a subset $\Omega$ of a Euclidean space $\mathbb{R}^n$ equipped with a metric to which the term ``shortest distance" refers, and where  the term ``shortest line" refers to a Euclidean line, or a Euclidean segment. One should remember however that the axioms for a metric space, as we intend them today, were formulated (by M. Fr\'echet) only  in 1907, that is, seven years after Hilbert stated his problem. From Hilbert's own elaboration on his problem \cite{Hilbert-Problems}, it appears that he  was thinking, rather than of a metric space, of a geometrical system, like Euclid's axiomatic system (or like his own axiomatic system); that is, a system consisting of undefined objects (points, lines, congruence, etc.) and axioms that make connections between these objects. In fact, Hilbert's problem can be included in a perspective that was dear to him on the axioms of geometry, and, more generally, on the foundations of mathematics.\footnote{We can recall here Hilbert's own words on a system of axioms for a geometry: ``It must be possible to replace in all geometric statements the words point, line, plane by table, chair, beer, mug", cf. Hilbert's \emph{On the axiomatic method in mathematics} (recalled in \cite{Ried} p. 57). One may also recall here that others, before Hilbert (and probably since Euclid's epoch), had the same point of view on the axiomatic method, even if they did not express it with the same words.} Indeed, in his comments about the problem, Hilbert gave the following formulation of the question\footnote{The translation is in \cite{Hilbert-Problems}.}:  
\begin{quote}\small
 In analogy with the way Lobachevsky (hyperbolic) geometry sits next to Euclidean geometry, by denying the axiom of parallels and retaining all the other axioms of Euclidean geometry, explore other suggestive viewpoints from where geometries may be devised which stand, from such a point of view, next to Euclidean geometry.
\end{quote}

This sentence clearly shows that Hilbert suggested to work out a theory in which one negates one of the axioms (other than the axiom of parallels) of Euclidean geometry, while keeping the other axioms untouched. Such an interpretation of Hilbert's Problem IV led to interesting developments which we shall mention in the next section, but of course no geometry whose importance is comparable to the Lobachevsky geometry was discovered by keeping, among the axioms of Euclidean geometry, the parallel axiom untouched and denying another one.

In fact, Hilbert gave a more precise formulation on that problem in the axiomatic framework, by adding that he asks for 
\begin{quote}\small
a geometry in which all the axioms of ordinary Euclidean geometry hold, and in particular all the congruence axioms, except the one of the congruence of triangles and in which, besides, the proposition that in every triangle the sum of two sides is greater than the third is assumed as a particular axiom.\footnote{We note that the word ``congruence", in Hilbert's point of view,  designates a relation between primary objects like segments, angles, triangles, etc. satisfying certain properties (some of which are axioms), and that this word does not necessarily refer to an equivalence relation in a metric sense, that is, an isometry. In fact, as in Euclid's axioms, there is no distance or length function involved in Hilbert's axioms. Congruence can be thought of as an undefined notion.} 
\end{quote}

Hilbert's problem, in this setting, can be formulated as follows: Consider the axiom system of Euclidean geometry, drop from it all the congruence axioms that contain the notion of angle\footnote{For instance, in Hilbert's axioms for Euclidean geometry, the following axiom is part of the set of congruence axioms for angles: 
\begin{quote}\small
 Given an angle  $\widehat{ABC}$ and given a point $B'$ and a ray $B'C'$ starting at $B'$, and given a choice of a side on the line $B'C'$, there exists  a ray $B'A'$ starting at $B'$, with $A'$ being on the chosen side of $B'C'$, such that  $\widehat{A'B'C'}\equiv \widehat{ABC}$. 
\end{quote}
Note that in Euclid's \emph{Elements}, this is a theorem (Proposition 23 of Book 1).} and replace them by an axiom saying that in any triangle the sums of the lengths of any two sides is not smaller than the length of the third side (that is, the triangle inequality). Then:
\begin{enumerate}
\item characterize the geometries satisfying these conditions;
\item study individually such geometries.
\end{enumerate}

We shall elaborate more on the axiomatic point of view in the next section.

 \section{Some early works on Hilbert's Problem IV}\label{s:early}
 
Among the major mathematicians who worked on Hilbert's Problem IV, there are Hilbert himself and his students Dehn\footnote{\label{n:Dehn} Max Dehn (1878-1952) is one of the main founders of combinatorial topology and combinatorial group theory. His thesis, defended in 1900, was in the setting of the foundations of geometry, and it is related to Hilbert's Problem IV, in its extended form discussed by Hilbert. He constructed a geometry in which there are infinitely many lines passing through a point and disjoint from another line and in which Legendre's theorem stating that the sum of the angles in any triangle is at most $\pi$ fails (therefore this geometry is different from hyperbolic geometry).  Dehn's geometry is non-Archimedean, and in this work Dehn constructed a \emph{non-Archimedean Pythagorean field}. The same year, Dehn, who was only 22 years old, solved Hilbert's Problem III which we already mentioned in Footnote \ref{Problems}. The problem asks whether two 3-dimensional (Euclidean) polyhedra that have the same volume are scissors-equivalent, that is, whether they can be obtained from each other by cutting along planes and re-assembling. The corresponding result in dimension 2 is true and is known as the \emph{Bolyai-Gerwien Theorem}. Dehn showed that in dimension 3 the answer is no, which was Hilbert's guess. In his work on that problem, Dehn introduced  a scissors-equivalence invariant for 3-dimensional polyhedra, which is now called the \emph{Dehn invariant}, and he showed that this invariant does not give the same value to the tetrahedron and the cube.  In 1965, J.-P. Sydler -- an amateur mathematician -- completed Dehn's theory by showing that two polyhedra are scissors-equivalent if and only if they have the same volume and the same Dehn invariant.}, Hamel\footnote{Georg Hamel (1877-1954) entered G\"ottingen University in 1900, and in 1901 he was awarded his doctorate, under the supervision of Hilbert. Hamel is probably best known for his contribution on the so-called \emph{Hamel basis}, a work published in 1905, in which he gave a construction of a basis of the real numbers as a vector space over the rationals. This construction is important because it involves one of the earliest applications of the axiom of choice. Hamel later made a name for himself in the field of Mechanics.} and Funk\footnote{Paul Funk (1886-1969) is mostly known for his work on the calculus of variation and for his introduction of the \emph{Funk transform} (also called the  \emph{Minkowski-Funk transform}),  an integral transform obtained by integrating a function on great circles of the spheres. The subject of his dissertation work (1911), surfaces with only closed geodesics, which he wrote under Hilbert, was further developed by Carath\'eodory. We owe him the Funk metric, an asymmetric metric defined on convex subsets of $\mathbb{R}^n$ which is a non-symmetric version of the Hilbert metric satisfying the requirements of Hilbert's Problem IV. This metric is studied in the chapter \cite{PT} of this Handbook.}. The title of Dehn's dissertation  is \emph{Die Legendre'schen S\"atze \"uber die Winkelsumme im Dreieck} (1900) (Legendre's theorems on the sum of the angles in a triangle), the title of the one of Hamel is \emph{\"Uber die Geometrieen in denen die Geraden die k\"urzesten sind} (1901) (On the geometries where the straight lines are the shortest), and the title of the one of Funk is \emph{\"Uber Fl\"achen mit lauter geschlossenen geod\"atischen Linien} (On surfaces with none but closed geodesic lines (1911)). All these works are related to Problem IV. Among the other early works that are related to Hilbert's problem, one also has to mention those of Finsler\footnote{Paul Finsler (1894-1970) wrote his dissertation in G\"ottingen in 1919 under the supervision of Carath\'eodory. The title of his thesis is \emph{\"Uber Kurven and Fl\"achen in allgemeinen R\"aumen} (On Curves and surfaces in general spaces). According to Busemann \cite{Busemann-F}, Finsler was not the first to study Finsler spaces. These spaces had been discovered by Riemann, who mentioned them in his lecture \emph{\"Uber die Hypothesen, welche der Grometrie zu Grunde liegen} (On the Hypothesis which lie at the foundations of Geometry) (1854) \cite{Riemann}. To support this claim, Busemann recalls that the goal that Riemann set for himself was the \emph{definition and discussion of the most general finite-dimensional space in which every curve has a length derived from an infinitesimal length or line element}. The name \emph{Finsler space} is due to E. Cartan, after his paper \emph{Les espaces de Finsler} \cite{Cartan1934} published in 1934. In their paper \cite{BP1993}, Busemann and Phadke write the following: ``Finsler was the first who investigated non-Riemannian spaces (under strong differentiability hypotheses) guided by Carath\'eodory, whose methods in the calculus of variations form the basis of Finsler's thesis, which has no relation to our $G$-spaces. The misnomer was caused by the fact that Finsler's thesis (1918) was inaccessible until 1951, when Birkh\"auser reissued it unchanged".} and of Berwald.\footnote{Ludwig Berwald (1883-1942) obtained his doctorate at the Royal Ludwig-Maximilians-University of Munich under the guidance of Auriel Voss. The title of his dissertation was: \emph{\"Uber die Kr\"ummungseigenschaften der Brennflachen eines geradlinigen Strahlsystems und der in ihm enthaltenen Regelfl\"achen} (On the curvature properties of the caustics of a
system of rays and the ruled surfaces contained in it). Berwald is considered as one of the founders of Finsler geometry.} We shall report on some of their works. 
 
Let us now say a few words on the work of Hilbert on the foundations on mathematics, in order to explain the axiomatic approach that he alludes to in the comments on his fourth problem. 

In the various revised editions of his \emph{Grundlagen der Geometrie}\footnote{The first edition appeared in 1899, that is, one year before Hilbert formulated his \emph{Problems}.} (Foundations of Geometry) \cite{Hilbert-Grund},  Hilbert worked out several sorts of geometries which he termed as ``non-Euclidean",\footnote{The term non-Euclidean is used here in a broad sense, and not only in its classical sense which usually includes only hyperbolic and spherical geometry.} a word he used in a very broad sense, in the formulation of his fourth problem. This includes first the well-known \emph{non-Archimedean} geometry,\footnote{Guiseppe Veronese (1854-1917)  constructed, in 1889 (that is, a few years before Hilbert did), a non-Archimedean geometry. Hilbert mentions the work of Veronese in the statement of his problem \cite{Hilbert-Problems}.} in which all the axioms of Euclidean geometry are kept untouched except the Archimedean axiom. We already mentioned that non-Archimedean geometry was also the subject of the doctoral thesis of Dehn (1902). Then, Hilbert constructed a \emph{non-Arguesian} geometry, which is a plane geometry in which the theorem of Desargues\footnote{\label{f:Desargues}Since this theorem plays a central role in the work of Hilbert and in the work of Busemann on axioms, we recall one version of it. Consider in the Euclidean plane (or  in the projective plane, for a more general version that is more convenient for the special cases where the lines are parallel) two triangles $abc$ and $ABC$. We say that they are \emph{in axial perspectivity} if the three intersection points of lines $ab\cap AB, ac\cap AC, bc\cap BC$ are on a common line. We say that the three triangles are \emph{in central perspectivity} if the three lines $Aa, Bb, Cc$ meet in a common point. Desargues theorem then says that for any two triangles, being in axial perspectivity is equivalent to being in central perspectivity.} of projective geometry fails,\footnote{A model of non-Arguesian geometry, described by Hilbert in the first editions of his  \emph{Grundlagen der Geometrie}, is obtained by taking an ellipse in the Euclidean plane and replacing the Euclidean segments intersecting it by arcs of circles that pass by a common fixed point $P$. Thus, in this geometry, a \emph{line} is taken in the usual sense if it does not intersect the (interior of the) ellipse, and in the case where it intersects it, it is made out of three pieces, two half-lines, outside the ellipse, connected by an arc of a circle inside the ellipse. The group of Euclidean axioms which Hilbert calls the ``projective axioms" are satisfied by such a geometry provided the point $P$ is not too close of the ellipse. In later editions, Hilbert replaced this model by a model found by F. R. Moulton, cf. \cite{Moulton}.} and a \emph{non-Pascalian} geometry, in which the theorem of Pascal fails.\footnote{Pascal's theorem in projective geometry says that given a hexagon in the projective plane, there is an equivalence between the following two statements: (1) the hexagon is inscribed in a conic; (2)  The interesections of the lines containing pairs of opposite sides are aligned.} To define these geometries, one has to construct new number fields, in particular a field of \emph{non-Archimedean numbers} and a field of \emph{non-Pascalian numbers}. A non-Pascalian geometry is also non-Archimedean. Hilbert also considered a geometry he called \emph{non-Legendrian}, in which it can be drawn from a point infinitely many parallels to a line that does not contain it, and where, like in spherical geometry, the angle sum in a triangle is greater than two right angles. He then considered a geometry he called \emph{semi-Euclidean}, in which one can draw from a point infinitely many parallels to a line not containing it, and where the angle sum in a triangle is equal to two right angles. Hilbert also considered a \emph{non-Pythagorean geometry}. In his paper \emph{\"Uber eine neue Begr\"undung der Bolyai-Lobatschefkyschen Geometrie} \cite{Hilbert-Bolyai}, he constructed a non-Archimedean Lobachevsky geometry. Another ``non-Euclidean" geometry that Hilbert mentioned is the Minkowski geometry, which he described as a geometry in which all the axioms of Euclidean geometry are satisfied except a ``triangle congruence" axiom. In this geometry, the theorem stating that the angles at the basis of an isosceles triangle are equal is not satisfied and the theorem stating that in any triangle the sum of the lengths of two sides is less than the length of the third side is taken as an axiom. In his comments on his fourth problem, Hilbert mentioned a geometry he described in his paper \cite{Hilbert1895} where the parallel axiom is not satisfied but all the other axioms of a Minkowski geometry are satisfied.

In this paraphernalia of geometries, only a few satisfy the requirements of Problem IV, but it is clear from the statement of the problem and the work that Hilbert did at the same period that he was thinking about new geometries in which a variety of other properties could replace the fact that ``shortest distances are realized by straight line".

Poincar\'e commented on these geometries and on the number fields to which they give rise in his review of Hilbert's \emph{Grundlagen der Geometrie}, \cite{Poincare-review-H}. He writes (p. 2):
 \begin{quote}\small
Many contemporary geometers [...] in recognizing the claims of the two new geometries [the hyperbolic and the spherical, which Poincar\'e calls the Lobachevsky and the Riemann geometries] feel doubtless that they have gone to the extreme limit of possible concessions. It is for this reason that they have conceived what they call \emph{general geometry}, which includes as special cases the three systems of Euclid, Lobachevsky, and Riemann, and does not include any other. And this term, \emph{general} indicates clearly that, in their minds, so other geometry is conceivable. 
\\
They will loose this illusion if they read the work of Professor Hilbert. In it they will find the barriers behind which they have wished to confine us broken down at every point.
\end{quote}

In the conclusion to his comments on Problem IV, Hilbert declares that it would be desirable to make a complete and systematic study of all geometries in which shortest distances are realized by straight lines. The conclusion of the problem in its French version is different from that of the English one: in the former, Hilbert includes the problem in the setting of the calculus of variations, and he refers to the work of Darboux that we already alluded to. The reader can find an English translation of the conclusion in the Appendix to the present paper.

It may be worth mentioning that in order to make the relation between the axiomatic abstract setting  and the metric setting, one may refer to Euclid. We already recalled that there is no distance function involved in Euclid's \emph{Elements}, but if we introduce this metric language, it is a consequence of Euclid's axioms (in particular Axiom 2) that on a line, distance is additive, and it follows from Proposition 20 of Book I that if between three points distance is additive then the three points are on a line (see \cite{Heath}). This implies that the Euclidean lines in the sense of the axioms of geometry are the (geodesic) lines in the metric sense. 

We now say a few words on other early works on Problem IV. 

Hamel considered in \cite{Hamel1903} metric spaces $(\Omega, d)$ where $\Omega$ is an open subset of a projective space $\mathbb{P}^n$ for some $n>1$ and where $d$ satisfies the following strong version of Hilbert's requirements (stated here in modern terms, following Busemann \cite{Busemann1970}):
\begin{enumerate}
\item \label{t1} $(\Omega, d)$ is a geodesic metric space and the projective lines are geodesics;
\item \label{t2} the closed balls are compact;
\item \label{t3} if $x,y,z$ do not lie on a projective line, then $d(x,y)+d(y,z)> d(x,z)$.
\end{enumerate}
Note that Conditions (\ref{t1}) and (\ref{t3}) imply that the metric $d$ is uniquely geodesic.

Under these hypotheses, Hamel proved that the set $\Omega$ is necessarily convex. More precisely, he  proved the following:
\begin{theorem} \label{th:Hamel}
If the metric space $(\Omega,d)$ satisfies (\ref{t1}) to (\ref{t3}), then  $\Omega$ is of one of the following two types:
\begin{enumerate}
\item \label{c1} $\Omega$ can be mapped homeomorphically onto the projective space $ \mathbb{P}^n$ in such a way that the maximal geodesics of $\Omega$ are mapped to great circles and all great circles have equal lengths;
\item \label{c2} $\Omega$ can be mapped homeomorphically onto an open convex subset of affine space $\mathbb{A}^n$ in such a way that the maximal geodesics of $\Omega$ have infinite length and are mapped to the intersection of $\Omega$ with affine lines in $\mathbb{A}^n$.
\end{enumerate}
\end{theorem}
In particular, closed geodesics and maximal open geodesics cannot coexist in $(\Omega,d)$ .

Let us note that the two geometries of constant curvature +1 and -1 are examples of metric spaces described in Theorem \ref{th:Hamel}, namely, spherical geometry is an example of type (\ref{c1}), and hyperbolic geometry, through the Klein-Beltrami model, is an example of type (\ref{c2}).

Hamel considered only smooth metrics, and his methods are variational. In his thesis \cite{Hamel1901} (1901), he considered non-necessarily symmetric metrics. In his later paper \cite{Hamel1903} (1903) he considered also non-symmetric metrics although he was not very explicit about that.  He refers to symmetric metrics by the terms ``Stark Monodromieaxiome"  and to non-symmetric metrics by ``Schwache Monodromieaxiome".\footnote{That is, ``strong monodromy hypothesis" and ``weak monodromy hypothesis" respectively. The reason for this terminology comes from dimension two, where Hamel (and Funk after him) use polar coordinates. In this case, the Lagrangian depends on $r$ and $\theta$. The periodicity of the angle variable $\theta$ imposes that its value is invariant by adding integer multiples of $2\pi$ (this is the weak monodromy hypothesis).  In the case of a symmetric metric, the value of $\theta$ must be 
invariant by adding integer multiple of $\pi$ (this is the strong monodromy hypothesis). I owe this remark to Marc Troyanov.}  The main result he obtained in that paper is a characterization of the differential of arclength, for a general smooth (symmetric and non-symmetric) metric  satisfying Hilbert's requirement, in dimensions 2 and 3. Hamel's result is presented in the chapter \cite{Tro-Finsler} in this Handbook.
 
 At its publication, Hamel's paper \cite{Hamel1903} was considered as a solution to Hilbert's problem, under the given restrictive assumptions. Concerning Hamel's method, in his 1974 report \cite{Busemann1974} on Hilbert's Problem IV, Busemann writes (p. 137):
 \begin{quote}\small  
 We do not reproduce it here, since it is not based on general ideas like integral geometry and the expression given by Hamel is not illuminating. Also, inevitably in view of the time it was written, the discussion of completeness does not satisfy our modern requirements; the necessary concepts did not yet exist. However, the immense number of possibilities becomes quite clear from \cite{Hamel1903}.
 \end{quote}

Busemann, in \cite{Busemann1955}, gave a proof of the same result, putting it in a modern setting and using completely different methods, which are purely metrical. He gave a non-symmetric version which he attributed to Hamel. A proof  of Hamel's result is contained in Busemann \cite{Busemann1970} p. 37, which Busemann states as follows in the setting of non-symmetric metrics:

\begin{theorem}
Under the hypotheses of Theorem \ref{th:Hamel} but without the assumption that the metric is symmetric and where Condition (\ref{t2}) is replaced by the requirement that the \emph{right} closed balls are compact,\footnote{For a non-symmetric metric $d$, one has to distinguish between the \emph{right open ball} of center $x$ and radius $r$, $B^+(x,r)= \{y\vert d(x,y<r\}$ and the \emph{left open ball} of center $x$ and radius $r$, $B^+(x,r)= \{y\vert d(y,x)<r\}$. The right and left closed balls are defined by replacing the strict inequalities by large inequalities.}  the space $(\Omega,d)$ is of one of the following two types:
\begin{itemize}
\item $\Omega= \mathbb{P}^n$, and each segment of a projective line in $\mathbb{P}^n$, traversed in either direction, is a geodesic. Furthermore, all such geodesic are projective lines, traversed in either direction, have the same length.
 \item $\Omega$ is a convex set and contained in some affine $n$-dimensional subspace $\mathbb{A}^n$ of $\mathbb{P}^n$. 
 \end{itemize}
 \end{theorem}
 
  Finally, we mention a result of Berwald related to Hilbert's problem IV. Berwald proved in  \cite{Berwald1929} (with a simpler proof given by Funk in \cite{Funk1935}) that 2-dimensional Finsler spaces of constant flag curvature\footnote{There are several notions of curvature in the setting of Finsler spaces, and an important one is the notion of flag curvature, which is a kind of Finsler analogue of the sectional curvature of Riemannian geometry. We refer the reader to \S 41 of Busemann's \emph{Geometry of Geodesics} \cite{Busemann1955} for a discussion of curvature in Finsler spaces and to the chapter \cite{Tro-Finsler} in this handbook.} that satisfy the requirements of Hilbert's Problem IV  are characterized by the property that any isometry of a geodesic onto another (or onto itself) is the restriction of a projectivity. This result generalizes Beltrami's result which he obtained in the setting of Riemannian geometry. A proof of this result is given in the chapter \cite{Tro-Finsler} of this Handbook.
  
  The works of Hamel, Funk and Berwald had a great influence on Busemann, who made a new definite breakthrough on the problem, which we describe in the next section.

  \section{Busemann's approach and Pogorelov's solution}
  
  We start with a reformulation of the problem.
  
We already mentioned that Hilbert's Problem IV, in its original formulation, is sometimes considered as too broad, and several mathematicians in order to work on it, made it more precise. Busemann, who spent a large part of his intellectual activity thinking and working on that problem reformulated it as follows (see \cite{Busemann1974}): 
\begin{quote}\small
The fourth problem concerns the geometries in which the ordinary lines, i.e. lines of an n-dimensional (real) projective space $\mathbb{P}^n$ or pieces of them, are the shortest curves or geodesics. Specifically, Hilbert asks for the construction of these metrics and the study of the individual geometries. It is clear from Hilbert's comments that he was not aware of the immense number of these metrics, so that the second part of the problem is not at all well posed and has inevitably been replaced by the investigation of special, or special classes of, interesting geometries. 
\end{quote}

  The second part of the problem is indeed very wide, and because of that it is likely to remain open forever. Busemann writes in \cite{Busemann-Review}:
  \begin{quote} \small
The discovery of the great variety of solutions showed that Part (2) of Problem 4 is not feasible. \emph{It is therefore no longer considered as part of the problem}. But many interesting special cases have been studied since 1929.
  \end{quote}
    Part (2) of the problem is nevertheless fascinating, in many cases of the individual geometries. There are various examples where new phenomena can be found, even in the following classes of classical metric spaces:
   \begin{enumerate}
   \item Minkowski metrics (that is, translation-invariant metrics on $\mathbb{R}^n$ defined by norms, which may be non-symmetric);\footnote{The Minkowski metrics are the translation-invariant metrics associated to finite-dimensional normed spaces, except that the norm function (and therefore the distance function) is not necessarily symmetric. See the chapter \cite{PT-Minkowski} in this Handbook.}
   \item Funk metrics;\footnote{For a survey on Funk geometries, see the chapter  \cite{PT} in this Handbook.}
   \item Hilbert metrics.
   \end{enumerate}

Each of the geometries for which the Euclidean straight lines are geodesics defines a new world, and one may try to answer in that world numerous questions regarding triangles, boundary structure, trigonometry, parallelism, perpendicularity, horocycles, compass and straightedge constructions, and there are many others. The two-and three-dimensional cases are particularly worth investigating in detail.

  According to Busemann \cite{Busemann-Review}, at the time where Hilbert proposed his problem, the only classes of metrics (``besides the elementary ones")  that were known to satisfy the requirements of that problem were the Minkowski metrics, introduced around 1890 by H. Minkowski for their use in number theory, and the Hilbert metrics. The latter were discovered by Hilbert in 1895 \cite{Hilbert1895}, as a generalization of the Klein-Beltrami model of hyperbolic geometry. We note that both classes of metrics are Finsler and almost never Riemannian.\footnote{The exceptions in both cases are represented by ellipses: the case where the convex set  on which the Hilbert metric is defined is an ellipse, and the case where the unit ball of the Minkowski metric is an ellipse centered at the origin.} Indeed, it is a consequence of Beltrami's work that we already mentioned, that the only  Riemannian metrics that satisfy Hilbert's requirements are the metrics of constant curvature. The Minkowski and the Hilbert metrics are mentioned in Hilbert's comments on Problem IV.  Busemann notes furthermore that Minkowski geometries satisfy the parallel axiom and that Hilbert geometries generalize the hyperbolic one.\footnote{Recall that the Hilbert metric of the disc is the hyperbolic metric.}  Busemann adds that probably Hilbert did not think of ``mixed situations", which Busemann himself considered in his various books and papers and which we recall now.
  
  Busemann introduced new classes of metrics which satisfy Hilbert's requirements and which are combinations of Minkowski and Hilbert metrics. In his book \cite{Busemann1955} (p. 111 ff.), he first makes the simple remark that the metric $L(p,q)+ e(p,q)$ where $L$ is the Hilbert metric and $e$ is the Euclidean metric has the required properties, and then he defines several variations on this metric. In particular, he considers a combination of the Euclidean metric with the Hilbert metric restricted to the intersection of the convex set with a flat subspace.
 He then defines (in dimension 2) several classes of combinations which depend on parameters and which satisfy interesting parallelism properties (in the sense of hyperbolic geometry). In \cite{Busemann1974}, returning to these examples, he considers the explicit metric $e(x,y) + \vert e^{x_{n}}-e^{y_{n}}\vert $ on affine space with affine coordinates $x=(x_1,\ldots x_n)$, and he notes that there are variations on this example where instead of the exponential function one can take any monotone function. He then adds the following:
 \begin{quote} \small
Many such functions can be added, and watching convergence even infinitely many. Also, the addition can be replaced by an integration over a continuous family. Of course, it is very easy to generate nonsymmetric distances by this method, for example by modifying nonsymmetric Minkowski distances. It seems clear that Hilbert was not aware of these simple possibilities. 
\end{quote}

For arbitray open convex subsets of $A^n$, Busemann also exhibits ``infinitely many essentially different Desarguesian metrics", namely, he takes the metric $e(x,y)+ \vert \tan x_n-\tan y_n\vert$ defined on the strip $-\pi/2 <x_n < \pi/2$, and again, he notes that the only property used is that $\tan t$ increases monotonically from $-\infty$ to $+\infty$ in $(-\pi/2, \pi/2)$.
  
 One of the major achievements of Busemann was to put Hilbert's Problem IV in an adequate metric setting. In particular, he introduced the terminology of a \emph{$G$-space} and of a  \emph{Desarguesian space} which we recall now, and over several decades he did extensive work on the geometry of these spaces.  At several places, he gave slightly different but overall equivalent definitions of $G$-spaces and of Desarguesian.

In  \cite{Busemann1942} and \cite{Busemann1943}, Busemann gave the following definition, see also Busemann's book \emph{The geometry of geodesics} \cite{Busemann1955}. 
\begin{definition}[$G$-space]A $G$-space\footnote{The expression $G$-space stands for ``geodesic space". The concept was introduced quite early by Busemann. In their paper \cite{BP1993}, p. 181, Busemann and Petty recall that ``The $G$-spaces were first introduced without this name by Busemann in his thesis in 1931." This notion turned out to be of paramount importance, and a lot of work has been done, by Busemann and others, especially around the conjecture of Busemann saying that every finite-dimensional $G$-space is a topological manifold. Busemann himself proved the conjecture in dimensions 1 and 2 (\cite{Busemann1955} \S 9 and \S 10), Krakus \cite{Krakus1968} proved it in dimension 3 using a theorem of Borsuk \cite{Borsuk} and P. Thurston \cite{Thurston1996} in dimension 4. It may useful to recall in this respect that $G$-spaces are metric and separable, and therefore the usual concept of dimension in the sense of Meyer-Urysohn applies to them, cf. \cite{HW}.} is a metric space $(X,d)$ satisfying the following properties:
\begin{enumerate}
\item (Finite compactness) Every bounded infinite set has an accumulation point.
\item (Menger convexity) For any $x$ and $z$ in $X$, there exists a point $y$ which is between $x$ and $z$, that is, $y$ different from $x$ and $z$ and satisfies $d(x,y)+d(y,z)=d(x,z)$.
\item (Local extendability) For every point $p$ of $X$ there exists an open ball $B$ centered at $p$ such that for  every $x$ and $y$ in $B$ there exists a point $z$ in $B$ such that $y$ is between $x$ and $z$. 
\item (Uniqueness of extension) If $y$ is between $x$ and $z_1$ and between $x$ and $z_2$ and if $d(y,z_1)=d(y,z_2)$ then $z_1=z_2$. 
\end{enumerate}
\end{definition}

In particular, in a $G$-space, geodesics are extendable, and the extension is unique.

 Busemann gave then the following definition of a  Desarguesian space (see \cite{Busemann1955} and \cite{Busemann-Review}):

\begin{definition}[Desarguesian space]
A Desarguesian space is a metric space $(R,d)$ satisfying the following properties:
\begin{enumerate}
\item $R$ is an open nonempty subset of projective space $\mathbb{P}^n$ equipped with a metric $\rho(x,y)$ whose associated topology is the one induced from its inclusion in $\mathbb{P}^n$.
\item The closed balls in $R$ are compact.
\item  Any two points in $R$ can be joined by a geodesic segment, that is, a curve of length $\rho(x,y)$. (Here, the length of a curve is defined as usual as the supremum over all subdivisions of the image set of a segment of the set of lengths of polygonal curves joining the vertices of the subdivision;  in modern terms, a metric where any two points can be joined by a geodesic segment is said to be \emph{geodesic}).
\item For any three points $x,y,z$ in $\mathbb{P}^n$ that are not collinear, the strict triangular inequality $\rho(x,y)+ \rho(y,z)< \rho(x,z)$ holds between them; that is, the three points do not belong to a geodesic.
\end{enumerate}
\end{definition}

\begin{remark}[Non-symmetric metrics]
In the later generalizations of $G$-spaces and of Desarguesian space to non-symmetric metrics that Busemann considered (see e.g. \cite{Busemann1970}), the condition that the closed balls are compact is replaced by the condition that the right closed balls are compact.
\end{remark}

Thus, whereas the definition of a $G$-space concerns general metric spaces, the one of a Desarguesian space concerns metrics on subsets of projective space. A Desarguesian space is a particular case of a $G$-space and it is also a particular case of a space satisfying the requirement of Hilbert's Problem IV. In a Desarguesian space, the geodesic joining two points is unique and is necessarity the affine geodesic. The terminology chosen by Busemann comes from the fact that in such a space,  Desargues' theorem\footnote{See Footnote \ref{f:Desargues}.} holds. Beltrami's theorem, which we already quoted, implies that in the Riemannian case, the only Desarguesian spaces are the spaces of constant curvature, that is, the Euclidean, hyperbolic or elliptic spaces.

The terminology used today for spaces satisfying Hilbert's problem requirements, especially in the context of Finsler geometry, is that of a \emph{projective metric}, or \emph{projectively flat metric} (cf. also Busemann in \cite{Busemann1953}).

 Let us note that Busemann was suspicious about the Finsler geometry viewpoint, in the sense in which this word was used in the early 1950s. His feeling about it is expressed in his survey paper \cite{Busemann-F} which starts as follows:
 \begin{quote} \small
 The term ``Finsler space" evokes in most mathematicians the picture of an impenetrable forest whose entire vegetation consists of tensors. The purpose of the present lecture is to show that the \emph{association of tensors} (or differential forms) \emph{with Finsler spaces is due to an historical accident, and that}, at least at the present time, \emph{the fruitful and relevant problems lie in a different direction}.
 \end{quote}
The paper ends with the following words:
 \begin{quote} \small
 This confirms that in spite of all the work on Finsler spaces we are now at a stage which corresponds to the very beginning in the development of ordinary differential geometry. Therefore the mathematician who likes special problems has the field. After sufficiently many special results have been accumulated someone will create the appropriate tools. At the present time it is difficult to guess what they will be beyond a vague feeling that some theory of integro-differential invariants will be essential.
 \end{quote}
          
Busemann interpreted Hamel's work in his setting of Desarguesian spaces; let us quote him from  \cite{Busemann1955} (p. 66):
  \begin{quote}\small  
  Whereas Hamel \cite{Hamel1903} has given a method for constructing all Desarguesian, sufficiently differentiable $G$-spaces, no entirely satisfactory infinitesimal characterization of these spaces, in terms of analogues to curvature tensors say, has ever been given. The freedom in the choice of a metric with given geodesics is for non-Riemannian metrics so great, that it may be doubted, whether there really exists a convincing characterization of all Desarguesian spaces. 
  \end{quote}
 
 Busemann also considered the axiomatic viewpoint. In \S 13 of his book \emph{The Geometry of Geodesics} \cite{Busemann1955}, he made relations between the metrical point of view on Problem IV and the \emph{foundations of geometry} point of view. The Section starts with the following introduction:
 \begin{quote} \small
 This section is closely related to the classical results of the \emph{foundations of geometry}, and the methods of this field are partly used here. We outline briefly the analogies as well as the differences of that work with the present.
 \end{quote}

 A. V. Pogorelov,\footnote{Pogorelov belonged to the Russian mathematical school, where he was taught by Alexandrov and Efimov. Before his work on Hilbert's Problem IV, Pogorelov had published several important works, including  solutions of some famous problems. In 1949, he gave a substantial generalization of Cauchy's rigidity theorem from convex polyhedra to convex surfaces. We recall that Cauchy's theorem (1813) states that any two convex polytopes in $\mathbb{R}^3$ which have congruent corresponding faces are congruent. Pogorelov's theorem states that any two closed isometric convex surfaces in $\mathbb{R}^3$ are congruent. The proof used a famous result of Alexandrov on gluing convex surfaces, which was used by him to prove that any surface homeomorphic to a sphere and equipped with a metric of positive curvature is isometric to a convex surface, 
 answering a famous problem which was asked by Weyl. The history of Cauchy's theorem is interesting; a mistake in Cauchy's proof was found by E. Steinitz in 1920, and corrected by him in 1928. Other improvements and extensions of the theorem are due to several people, among them M. Dehn (1916) and A. D. Alexandrov (1950).) In 1952, he published a solution of the so-called \emph{multi-dimensional regularity problem of the Minkowski problem}, stating that any convex surface whose Gaussian curvature is positive and is a $C^m$ function of the outer normal, for $m\geq 3$, is $C^{m+1}$-regular. (This problem has several later developments and generalizations, by Nirenberg, Yau and Cheng and others.)} in his monograph \cite{Pogorelov} that accounts for his solution of Hilbert's Problem IV, notes  that the Minkowski geometries (respectively the Hilbert geometries) realize the system of axioms of Euclidean geometry (respectively of Lobachevsky geometry) after we drop all the axioms of congruence involving the concept of angle and adding to them the triangle inequality axiom (which is one of Hilbert's requirements at the level of axioms). 
  He adds (without further elaboration) that Hilbert's problem can be formulated in a similar way for elliptic geometries. \'Alvarez Paiva in \cite{Alvarez2003} makes further comments on the elliptic geometry case.

    The idea of Pogorelov's proof of Hilbert's Problem IV came from Busemann, who, at the Moscow ICM (1966), gave a  construction of a Desarguesian space, which introduced ideas from integral geometry in the approach to Hilbert's problem. The method provided a metric on any bounded open convex set $B$ in a projective space $\mathbb{P}^n$. We already noted that (by a result of Hamel) a solution of Hilbert's problem is necessarily a metric defined either on $\mathbb{P}^n$ or on an open convex subset of affine space (sitting as the complement of a hyperplane in $\mathbb{P}^n$). We briefly recall Busemann's definition in case where $B=\mathbb{P}^n$. 
    
    Let $\mathcal{H}$ be the set of hyperplanes in $\mathbb{P}^n$. For any subset $X$ of $\mathbb{P}^n$, we let $\pi X$ be the set of hyperplanes that intersect $X$.    
    On the set $\mathcal{H}$, we define a nonnegative measure $m$ satisfying the following properties:
    
    \begin{enumerate}
    \item For any point $p$ in $\mathbb{P}^n$, $m(\pi \{p\})=0$.
   
    \item For every  nonempty open subset $X$ of $\mathbb{P}^n$, $m (\pi X)>0$.
    \item $m(\pi {\mathbb{P}^n})= 2k<\infty$.
    \end{enumerate}
    
    Then we notice that for every line $L$ in $\mathbb{P}^n$,  $\pi L= \pi \mathbb{P}^n$, which shows that $m (\pi L)=2k$.

    Now for every pair of distinct points $x$ and $y$ in $\mathbb{P}^n$, consider the line $L$ that contains them. The two points divide $L$ into two arcs, $A$ and $B$, and $m(\pi A) + m(\pi B)=
    m(\pi \mathbb{P}^n)=2k$. Therefore, one of the two values $m(\pi A)$ or $m(\pi B)$, say   $m(\pi A)$, is $\leq k$. Setting $\rho(x,y)=  m(\pi A)$ defines a Desarguesian metric on $\mathbb{P}^n$. 
    
    Explicit formulae can be obtained for such measures since the set of hyperplanes in projective space can easily be parametrized.
    
    One must remember here that Busemann, at the time he was invited to the Moscow 1966 ICM, was working in  relative isolation, and that the subject of metric geometry was not fashionable in the West. On the contrary, in the Soviet Union, the subject, which flourished independently, was considered as important, and it was represented by a host of very good and influential mathematicians including A. D. Alexandrov, V. A. Zalgaller, V. A. Toponogov, N. V. Efimov, Yu. G. Rechetnjak,  and A. V. Pogorelov himself. The interaction between the two schools, before Busemann's visit to Moscow, was poor.\footnote{We recall by the way that  in 1985, Busemann received the Lobachevsky prize, which was certainly the most prestigious reward in geometry that was given the Soviet Union (and, before that, in Russia). The prize was given to him ``for his innovative book 
\emph{The Geometry of Geodesics}" which he had written 30 years before. We also recall that the first recipients of this prize were Sophus Lie in 1897, Wilhelm Killing in 1900 and David Hilbert in 1903). The prize was awarded to Alexandrov in 1951.} 
 In last paper \cite{BP1993}, written in collaboration with Phadke and published in 1993, Busemann recalls his beginning in metric geometry. He writes the following (p. 181):
 \begin{quote}\small
 Busemann has read the beginning of Minkowski's \emph{Geometrie der Zahlen} in 1926 which convinced him of the importance of non-Riemannian metrics. At the same time he heard a course on point set topology and learned Fr\'echet's concept of metric spaces. The older generation ridiculed the idea of using these spaces as a way to obtain results of higher differential geometry. But it turned out that a few simple axioms on distance suffice to obtain many non-trivial results of Riemannian geometry and, in addition, many which are quite inaccessible to the classical methods.
 \end{quote}
 In the concluding remarks of the same paper, Busemann writes:
   \begin{quote}\small
   The acceptance of our theory by others was slow in coming, but A. D. Alexandrov and A. M. Gleason were among the first who evinced interest and appreciated the results.
   \end{quote}
   Let us now quote Pogorelov, from the introduction to his monograph \cite{Pogorelov}:  \begin{quote} \small The occasion for the present investigation is a remarkable idea due to Herbert Busemann, which I learned from his report to the International Congress of Mathematicians at Moscow in 1966. Busemann gave an extremely simple and very general method of constructing Desarguesian metrics by using a nonnegative completely additive set function on the set of planes and defining the length of a segment as the value of the function of the set of planes intersecting the segment.
  
  I suspected that all continuous Desarguesian metrics could be obtained by this method. The proof of this in the 2-dimensional case strengthened my belief in this conjecture and I announced a general theorem in \cite{Pogorelov1}. However it turned out later, on making a detailed investigation of the three-dimensional case, that the completely additive set function figuring in Busemann's construction may not satisfy the non-negativity condition. Therefore, the result given here, while preserving its original form, assumes that other conditions are satisfied.
  \end{quote}  
  
Busemann's construction mentioned by Pogorelov is based on the integral formula for distances defined as a measure on the set of planes that we mentioned above. The formula is contained in  earlier work of Blaschke \cite{Blaschke1936}, based on an integral formula which is due to Crofton\footnote{A formula for a metric on a subset of Euclidean space that uses a measure on the set of Euclidean hyperplanes, where the distance between two points $x$ and $y$ is equal to the measure of the set of Euclidean hyperplanes that meet the segment $[x,y]$, is sometimes referred to as a Crofton (or sometimes Cauchy-Crofton) formula. The name is after Morgan Crofton (1826-1915) who proved a result of integral geometry relating the length of a curve to the expected number of times a random line intersects it, see \cite{Santalo} p. 12-13.}  in the setting of geometric probability. We shall say a few more words on this formula below, since it was at the basis of most of the developments that followed.

   As already mentioned, Pogorelov's solution of Hilbert's Problem IV consisted in showing that \emph{every} Desarguesian metric is  given by Busemann's construction, but with a slight modification. Pogorelov obtained the result, first in dimension $n=2$ for a general continuous metric, and then for $n=3$, with a smoothness assumption. He also proved that in dimension 3, every metric that is a solution of Hilbert's Problem IV is a limit (in the sense of uniform convergence on compact sets) of a sequence of such metrics that are of class $C^\infty$. It turned out that the proof for $n=3$ needs more work than the one for $n=2$; in particular, a new definition of the measure $m$, allowing it to take negative values (but with still positive values on every set of the form $\pi T$, where $T$ is a segment) is needed. In fact, Pogorelov's solution (even in dimension 2) is variational and it uses the smoothness condition. For dimension 2, Pogorelov obtained the result for continuous metrics by using an approximation argument (which does not apply in dimensions $>2$).
       
    One should also note that the fact that in dimension $>2$ the measure $m$ can take negative values is also mentioned in Busemann's earlier paper \cite{Busemann1960} in which he showed that all Minkowski metrics which are smooth enough can be obtained by this construction.

        Pogorelov's exposition in the book \cite{Pogorelov} is elementary and pleasant. Busemann wrote a review on that book in \cite{Busemann-Review} in which he confirmed that Pogorelov's method of proof for $n=3$ works, with more technicalities, for $n>3$. It is fair to recall here that Hilbert, in the spirit of his time, restricted himself to dimensions 2 and 3.

 In their paper \cite{BP1984}, Busemann and Phadke work out a general and global\footnote{Let us recall that Beltrami's problem is indeed local.} version of  the 1865 Beltrami Theorem \cite{Beltrami1865} in the setting of Desarguesian spaces. They define 
 a general class of metric spaces which they call chord spaces. These are metric spaces possessing a distinguished class of extremals (that is, locally isometric maps from an interval into the metric space), called chords, which behave locally like the affine lines, but which are not the unique extremals of the metric. Chord spaces are the subject of the monograph \cite{BP1987}. 
Busemann and Phadke obtain in \cite{BP1984} the following results, which they ``believe to be the most general meaningful version of that theorem":
 \begin{theorem} A locally Desarguesian simply connected chord space is either defined in all of $S^n$ or is an arbitrary open convex set of an open hemisphere of $S^n$ (considered as the projective space $A^n$).
 \end{theorem}
  \begin{theorem} A simply connected locally desarguesian and locally symmetric $G$-space is Minkowskian, hyperbolic or spherical.
   \end{theorem}  

At the end of their paper \cite{BP1984}, Busemann and Phadke declare the following:
\begin{quote}\small
Surprisingly, one can clearly deduce from a paragraph (which is too long to be produced here in total) in Beltrami's paper \cite{Beltrami1868} that he would have welcomed our generalization of his theorem. He says that a theorem which holds under weaker hypotheses than stated has not been fully understood and that the proper generalization may involve the disappearance of some of the original concepts (in our case the Riemannian metric).
\end{quote}

            \section{Other works and other solutions}           
                 More general versions of Pogorelov's solution in the 2-dimensional case were later on obtained by completely different methods by Ambartzumian \cite{A1976} and by Alexander \cite{A1978}. Ambartzumian's method can be used only in dimension 2 and it gives the result of Pogorelov directly for continuous metrics, without passing by smooth metrics and without using any approximation argument. The method has a combinatorial character, and it uses the scissors-congruence of triangles which is the subject of Hilbert's Problem III which we already mentioned in Footnote \ref{Problems}. The fact that Ambartzumian's method does not work in dimension $>2$ can be compared to the fact that the solution to Hilbert's Problem III says that in dimension $>2$, two polyhedra with the same volume are not necessarily scissors-equivalent.

   Szab\'o, in a paper  \cite{Szabo} containing several new ideas on the subject, gave two new proofs of Pogorelov's theorem in a Finslerian setting, valid in all dimensions, one of them being elementary and based on ideas originating in Ambartzumian's work. Szab\'o's result is given in terms of a partial differential equation satisfied by the indicatrix of the Finsler structure, which in dimension $n$ is assumed to be of class $C^{n+2}$. The proof uses Fourier transform theory, and precisely, a relation between the problem of characterizing (symmetric) projective metrics of low regularity and that of characterizing the (distributional) Fourier transform of norms on finite-dimensional vector spaces.  For the relation between the Fourier transform and integral geometry on normed and projective Finsler spaces, see the paper \cite{AF} by \'Alvarez Paiva and Fernandes.

       A new point of view on Hilbert's problem, based on an integral formula in dimension 2, is given in the paper \cite{A1978} by Alexander, whose result can be considered as a general solution which does not even assume the continuity of the metric. In that paper, the author considers an indefinite metric (that is, a metric which does not necessarily separates points) on a 2-dimensional space homeomorphic to a Euclidean plane, equipped with a set of ``lines" (in an abstract sense),  which are homeomorphic to the real line. He shows that under very mild conditions (the metric is not assumed to be continuous and the set of lines are only assumed to satisfy Desargues theorem, there exists a unique Borel measure on the set of lines such that the distance between two points is given by integrating this measure. The methods are based on the author's combinatorial interpretation of the Crofton formula from integral geometry. This works includes Busemann's integral construction in the abstract setting of the axioms of geometry which was one of Hilbert's favorite settings. 
       
       The relation of Hilbert's Problem IV with the foundations of geometry is highlighted again in the paper \cite{Alexander} by Alexander in which this author extends his previous result to higher dimensions, using the notion of zonoids from convexity theory.\footnote{One first defines a \emph{zonotope} as a finite Minkowski sum of line segments. Equivalently, a zonotope is the affine projection of a high-dimensional cube. This can be seen as a generalization of the notion of polytope. A \emph{zonoid} is then a limit of zonotopes in the sense of Blaschke-Hausdorff convergence. Equivalently, a zonoid is the limit of sums of segments. We refer the reader to the survey \cite{Bolker} by E. D. Bolker.} Alexander works in the settings of \emph{hypermetrics}, that is, metrics $d$ such that for any points $P_1,\ldots,P_m$ in $\mathbb{R}^n$ and for any integers $N_1,\ldots,N_m$ satisfying $N_1+\ldots+N_m=1$, one has $\sum_{i<j} d(P_i,P_j)N_iN_j\leq 0$ (see \cite{Kelly}). This is a strong version of the triangle inequality (and for $m=3$, the property boils down to the triangle inequality). Alexander notes in \cite{Alexander} that in dimension 2, all solutions of Hilbert's problem are hypermetrics, and the main result of this paper is that there is a natural linear isomorphism between the cone of hypermetrics on $\mathbb{R}^n$ that are additive along lines (that is, hypermetrics that satisfy the requirement of Hilbert's Problem IV) and the cone of positive Borel measures $m$ on the set of hypersurfaces of $\mathbb{R}^n$ satisfying $m(\pi P)=0$ and $0<m(\pi [P,Q])=0<\infty$ for any $P$ and $Q$ in $\mathbb{R}^n$. Here, as before, for $X\subset\mathbb{R}^n$, , $\pi X$, denotes the set of hyperplanes that intersect $X$. The isomorphism is given  by a generalized Crofton-type formula $d(P,Q)=m(\pi [P,Q])$.
   Alexander's methods using zonoids were further generalized by Schneider (see \cite{Schneider2001} and \cite{Schneider2006}), with a result valid for a higher-dimensional analogue of Hilbert's problem involving area, volume, etc. and not only length. In the last section of this paper, we recall Busemann's formulation of this higher-dimensional analogue.

            In the paper \cite{AO1998}, Ambartzumian and Oganian  give (Theorem 1) a necessary and sufficient condition for a Finsler metric in dimension 2 to be a solution of Hilbert's problem.  In the same paper, they study a parametric version of the problem, where the convex domain which is the support of the metric (which is assumed to be centrally symmetric) depends continuously on parameters. They also consider the question of when the family of Finsler metrics obtained is a Minkowski metric. They obtain a complete solution in the case of a dependence on a single (1-dimensional) parameter. The smoothness condition is essential in their work.

Among the other recent works on the Finsler smooth case, we mention the works of Z. Shen and of R. L. Bryant. We recall that Funk, in his papers \cite{Funk} and \cite{Funk1935} gave a classification of Finsler metrics  on convex domains in the plane which have constant flag curvature.
In his paper \cite{Shen2003}, Shen described a class of smooth projective Finsler metrics of constant flag curvature using algebraic equations.
In the paper \cite{Funk1963}, Funk proved that the standard metric of the sphere $S^2$ is the only flat projective Finsler metric, satisfying some additional conditions. Bryant, in \cite{Bryant1997} and \cite{Bryant2002} generalized the result of Funk; he described a 2-parameter family of projective Finsler metrics on $S^2$ and he showed that this family gives all such metrics. In his work, Bryant made use of the symplectic geometry of the space of geodesics (the quotient of the unit tangent bundle by the geodesic flow).

     Finally, we mention that \'Alvarez Paiva gave in \cite{Alvarez2005} another point of view on Hilbert's Problem IV in the framework of symplectic and contact topology. He obtained a characterization of Finsler metrics which are projective (that is, which are solutions of Hilbert's problem) in terms of a class of symplectic forms on the space of lines in $\mathbb{R}^n$. At the same time, the result of \'Alvarez Paiva makes a relation between Hilbert's problem and the ``black and white" theory of Guillemin and Sternberg \cite{GS}. The results of \'Alvarez Paiva in his paper \cite{Alvarez2005}  were first announced in the paper \cite{AGS} by
 himself, Gelfand and Smirnov, 
     
 The paper \cite{Crampin2011} by M. Crampin contains reformulations and new versions of the works of Hamel and of \'Alvarez Paiva.
\section{Further developments and perspectives}
   
A problem like Hilbert's Problem IV never dies, and as we tried to show, it admits several generalizations and formulations in various settings. We mention that there are classical works on problems that are more general than Hilbert's problem, namely, problems where the metrics on subsets of the plane are to be found,  with prescribed geodesics which are curves that are not necessarily straight lines, and which satisfy certain conditions, e.g. they foliate the plane, and each of them homeomorphic to a real line that separates the plane into two components.  Blaschke and Bol consider such problems in \cite{Blaschke-Bol}. Before that, Darboux had already considered such a general problem, in the setting of the calculus of variations, see \cite{Darb} and he obtained local  existence results using theorems on partial differential equations.  Busemann and Salzmann also considered such general problems in \cite{Busemann-Salzman}.

 The case of nonsymmetric metrics is still wide open although, talking about this case,  Busemann mentions in \cite{Busemann1970} (p. 36): ``most contributors to Hilbert's Problem IV have treated this general case". In fact, Hilbert, when he formulated his problem, certainly had this case in mind, since he knew about the interesting examples of non-symmetric metrics discovered by Minkowski, and he was aware of the fact that extremal metrics obtained by the methods of the calculus of variation may lead to non-symmetric length functions. 
We already mentioned that Hamel, who proposed in his thesis \cite{Hamel1901} the first proof of Hilbert's problem for dimensions 2 and 3, considered the case of nonsymmetric metrics.  Szab\'o, in his paper  \cite{Szabo}, considered only the symmetric case and in the final section of that paper he outlined a sequel to that work which would treat the non-symmetric case; the sequel never appeared.
 Alexander, in the paper \cite{A1978} that we mentioned above discusses the nonsymmetric case but without being able to settle it. 

 Busemann also asked for the generalization of Hilbert's problem (see \cite{Busemann1960} and \cite{Busemann1961}) to the study of $k$-dimensional area in $n$-dimensional affine spaces for which the affine $k$-flats minimize area. We quote Busemann, from the introduction to \cite{Busemann1961}:
 \begin{quote}\small
 In his Problem IV, Hilbert [Nachr. Ges. Wiss. G\"ottingen Math.-Phys. Kl. 1900, 253-297] suggested studying all geometries in which the straight lines are the shortest connections or minimize length. Although the problem in this form proved too general to be fruitful, it has given rise to many investigations dealing either with special geometries or with special questions. The higher-dimensional analogue to Problem IV concerns $a$-dimensional metrics in subsets of the $n$-dimensional projective space in which the $a$-flats minimize area. Such geometries apparently have never been treated, at least in the large. The present paper represents a first step in this direction.
 \end{quote}

In the same paper, Busemann asks for the the study of $k$-dimensional areas that are given by a Crofton-type formula, similar to the one he gave for distances. We already mentioned advances on this problem, cf. the papers \cite{AF} by \'Alvarez Paiva and Fernandes, in which the authors prove the Crofton formulas for smooth projective metrics, and \cite{Schneider2001} by Schneider, in which the author establishes again the Crofton-type formula for the Holmes-Thompson area for Finsler projective hypermetrics.

     Besides the two questions in Hilbert's Problem, viz. the characterization of the metrics with the required properties and the study of such metrics individually, there is a third question, that lies in between, which is the question of finding general properties of such metrics. Busemann considered this question in various papers and books, and we mention as an example the following embedding result (part of it obtained in \cite{Busemann1955} p. 81, and it was completed in \cite{Busemann1970} p. 32): 
     \begin{theorem} For any $n$-dimensional Desarguesian $G$-space $R$ there is an $(n+1)$-dimensional Desarguesian space $R^*$ such that $R$ is a hyperplane in $R^*$ and the restriction of $R^*$ to $R$ is the given metric on $R$. 
     \end{theorem}
     
   Another general result  on Desarguesian spaces, also due to Busemann, is the following characterization of Hilbert and Minkowski geometries among all Desarguesian spaces (\cite{Busemann1970} p. 38):

   \begin{theorem} Let $(\Omega,d)$ be a non-compact Desarguesian space whose distance function $d$ is non-necessarily symmetric.
  Assume that any isometry between two lines in $\Omega$ is a projectivity. Then $(\Omega,d)$ is of one of the following two sorts:
  \begin{itemize}
  \item $\Omega$ is a convex subset of projective space and $d$ is its Hilbert metric;
  \item  $\Omega$ is affine space and $d$ is a Minkowski metric. 
  \end{itemize}
   \end{theorem}
 This is a special case of a theorem that Berwald obtained under an additional differentiability hypothesis and which we mentioned at the end of Section \ref{s:early}.

   The following result of Busemann gives a characterization of Funk and Minkowski geometries among Desarguesian spaces \cite{Busemann-homothetic}:
   
 \begin{theorem}
 Any Desarguesian space in which all the right spheres of positive radius around any point are homothetic is either a Funk space or a Minkowski space. 
   \end{theorem}

   We also mention a result by Busemann which is more in the axiomatic spirit of Hilbert. Busemann develops in \cite{Busemann1955} \S 23 and 24 a theory of parallelism between geodesics in metric spaces which generalizes the Euclidean theory.  He introduces in particular a \emph{parallel axiom} (see  \cite{Busemann1955} p. 141) and he proves the following (Theorem 24.1 of \cite{Busemann1955}):
   \begin{theorem}
A Desarguesian space in which the parallel axiom holds and where the spheres are convex is Minkowskian.
\end{theorem}

  Several of Busemann's students and collaborators worked on general properties of Desarguesian spaces. We mention in particular the following result of E. M.  Zaustinsky (one of Busemann's students) \cite{Zaustinsky1959}
  \begin{theorem}
A Desarguesian metrics defined on affine space and whose right and left open balls are all compact is  necessarily a Minkowski metric. 
\end{theorem}
B. B. Phadke \cite{Phadke1972} (who was also a student of Busemann),  proved the following.
  
  \begin{theorem}
  A Desarguesian space in which the equidistant loci to hyperplanes are hyperplanes is necessarily  a Minkowski space.
  \end{theorem}

Finally, we mention  the following:

\medskip

\noindent{\bf Problems}

\begin{enumerate}
\item \label{F2}  Modify Hilbert's fourth problem so that geometric circles replace lines, in spaces of constant curvature and on more general metric spaces. 
\item \label{F3}  Generalize Hilbert's Problem IV to the setting of metric spaces or convex subsets of metric spaces other than Euclidean spaces; that is, study the metrics on subspaces such that the geodesics of the ambient metric are geodesics for the new metric. (We recall that Busemann developed a theory of convexity in general metric spaces, see \cite{Busemann1955} and the exposition in \cite{P2}.)
\item \label{F4}  Study in particular the preceding problem on Teichm\"uller space equipped with its various metrics. 
\item \label{F5}  Study a complex analogue of Hilbert's problem.
\item \label{F6}  In the paper \cite{PB2006}, 
the authors used the great variety of explicit constructions of projective metrics to show that totally geodesic submanifolds (in that case, planes) are not minimal surfaces for the Hausdorff measure.  A nice question mentioned in that paper is to see whether the following holds: Suppose that for a given projective metric, planes (or hyperplanes) are minimal for the Hausdorff measure; is the metric then either flat or Riemannian ?
\item \label{F7}  The set of projective pseudo-metrics (pseudo distance functions) on $\mathbb{R}^n$ is a closed convex cone. Determine the extreme rays of this cone, for $n \geq 3$. 
\end{enumerate}

Problem ({F2}) was solved by \'Alvarez Paiva and Berck in \cite{PB2010}, for Finsler metrics and in the case of circles on the two-sphere. 

Regarding Problem (\ref{F3}),  we mention the work in \cite{PY1} and the survey in the chapter \cite{PY} in this Handbook, which are related to a version of Hilbert's Problem IV in the setting of constant curvature geometries. We study there two metrics (namely a Funk-type metric and a Hilbert-type metric) on open subsets of the sphere and of hyperbolic space such that the geodesics for the restriction of the ambient metric are geodesics for our two metrics. We can also mention the so-called Apollonian weak metric, which is a metric defined on open subsets of $\mathbb{R}^n$. In the case where the convex set is the open unit disc or the upper half-plane (considered as the ambient spaces of the Poincar\'e models of hyperbolic geometry), the non-Euclidean geodesics are also geodesics of the Apollonian metric, see \cite{PTA}, and the recent paper \cite{Yamada-A}.

The reason of the mention of Teichm\"uller space in Problem (\ref{F4}) is the amount of interest and of profound work that was done around it in the last six or seven decades. We note that there is a so-called \emph{Weil-Petersson-Funk metric} metric on Teichm\"uller space, which was first introduced by Yamada in  \cite{Yamada}, defined in analogy with the Funk metric on convex subsets of $\mathbb{R}^n$, (see the chapter \cite{MOY} in this Handbook and see also the survey \cite{Yamada-H}). Yamada also noticed that the Weil-Petersson geodesics are also  Weil-Petersson-Funk geodesics under a condition which concerns the monotonicity of the distance between the  Weil-Petersson-Funk geodesics and the relevant strata, parameterized by the  Weil-Petersson arclength $t$.  In fact, the distance function is Weil-Petersson convex in $t$, and the condition is about the function not having a critical point while traversing from one endpoint to the other. 
The conclusion is that except in some pathological situations, the Weil-Petersson-Funk metric is mostly projective, namely the geodesics of the underling geometry are also geodesics for this metric.

 There are several ways to formulate precisely the problem in (\ref{F5}), and we propose the use of complex geodesics (holomorphic maps from the  unit disc to a complex manifold). The Kobayashi metric could play the role played by the Euclidean metric in the real case.

Problem (\ref{F7}) is due to \'Alvarez Paiva. It is solved for $n=2$ in \cite{APB4}

 \'Alvarez Paiva studied the problem of extending Hilbert's problem to rank one symmetric spaces.  In a recent preprint,  \cite{Alvarez}, he shows that if $M$ is a manifold diffeomorphic to a compact rank-one symmetric space equipped with a reversible (i.e. symmetric) Finsler metric $F$ all of whose geodesics are closed and of the same length and if $L$ is a Finsler metric on $M$ whose geodesics coincide as unparametrized geodesics with those of $F$, then $L$ is the sum of a reversible Finsler metric and an exact 1-form.

Finally, \'Alvarez Paiva and Berck obtained recently a complete characterization of projective Finsler (not necessarily symmetric) metrics. The preprint \cite{Alvarez-B}  contains the details of the result for metrics defined on the whole projective $n$-space.
 
The story of Hilbert's Problem IV, of its various reformulations and of all the work done around it is an illustration of how mathematical problems can grow, transform and lead to new problems. It also shows how mathematicians are persistent men.

\section*{Appendix I: Hilbert's Problem 4}   
 
 (The translation is from the article in the 1902 Bulletin of the AMS.)
 
\noindent{\bf Problem 4: Problem of the straight line as the shortest distance between two points.}

\begin{quote}\small 
Another problem relating to the foundations of geometry is this: If from among the axioms necessary to establish ordinary Euclidean geometry, we exclude the axiom of parallels, or assume it as not satisfied, but retain all the other axioms, we obtain, as is well known, the geometry of Lobachevsky (hyperbolic geometry). We may therefore say that this is a geometry standing next to Euclidean geometry. If we require further that that axiom be not satisfied whereby, of three points on a straight line, one and only one lies between the other two, we obtain Riemann's (elliptic) geometry, so that this geometry appears to be the next after Lobachevsky's. If we wish to carry out a similar investigation with respect to the axiom of Archimedes, we must look upon this as not satisfied, and we arrive thereby at the non-Archimedean geometries which have been investigated by Veronese and myself. The more general question now arises: Whether from other suggestive standpoints geometries may not be devised which, with equal right, stand next to Euclidean geometry. Here I should like to direct your attention to a theorem which has, indeed, been employed by many authors as a definition of a straight line, viz., that the straight line is the shortest distance between two points. The essential content of this statement reduces to the theorem of Euclid that in a triangle the sum of two sides is always  greater than the third side -- a theorem which, as is easily seen, deals solely with elementary concepts, \emph{i. e.}, with such as are derived directly from the axioms, and is therefore more accessible to logical investigation. Euclid proved this theorem, with the help of the theorem of the exterior angle, on the basis of the congruence theorems. Now it is readily shown that this theorem of Euclid cannot be proved solely on the basis of those congruence theorems which relate to the application of segments and angles, but that one of the theorems on the congruence of triangles is necessary. We are asking then, for a geometry in which all the axioms of ordinary Euclidean geometry hold, and in particular all the congruence axioms except the one of the congruence of triangles (or all except the theorem of the equality of the base angles in the isosceles triangle), and in which, besides, the proposition that in every triangle the sum of two sides is greater than the third is assumed as a particular axiom.

One finds that such a geometry really exists and is no other than that which Minkowski constructed in his book, Geometrie der Zahlen,\footnote{Leipzig, 1896.} and made the basis of his arithmetical investigations. Minkowski's is therefore also a geometry standing next to the ordinary Euclidean geometry; it is essentially characterized by the following stipulations:

1. The points which are at equal distances from a fixed point $O$ lie on a convex closed surface of the ordinary Euclidean space with $O$ as a center.

2. Two segments are said to be equal when one can be carried to the other by a translation of the ordinary Euclidean space.

In Minkowski's geometry the axiom of parallels also holds. By studying the theorem of the straight line as the shortest distance between two points, I arrived\footnote{\emph{Math. Annalen}, vol. 46, p. 91.} at a geometry in which the parallel axiom does not hold, while all other axioms of Minkowski's geometry are satisfied. The theorem of the straight line as the shortest distance between two points and the essentially equivalent theorem of Euclid about the sides of a triangle, play an important part not only in number theory but also in the theory of surfaces and in the calculus of variations. For this reason, and because I believe that the thorough investigation  of the conditions for the validity of this theorem will throw a new light upon the idea of distance, as well as upon other elementary ideas, \emph{e.g.}, upon the idea of the plane, and the possibility of its definition by means of the idea of the straight line, {\it the construction and systematic treatment of the geometries here possible seems to me desirable}.
 
\end{quote}

\section*{Appendix II: Extract from the French version} 

{\bf Sur les probl\`emes futurs des math\'ematiques, 1902}

The conclusion of the French version is different from the English one:

\begin{quote}\small
In the plane, and assuming the continuity axiom, the problem in hand leads to a question treated by M. Darboux\footnote{\emph{Le\c cons sur la th\'eorie g\'en\'erale des surfaces, t. III, p. 54, 1894.}}: To determine all the calculus of variation problems in the plane whose solutions are all the plane straight lines; this question seems to me worth of and interesting generalizations.\footnote{Compare with the interesting research of M. A. Hirsch, (\emph{Math. Annalen}, t. XLIX et L).}

[Dans le plan et en admettant l'axiome de la continuit\'e, le probl\`eme dont il s'agit conduit \`a la question trait\'ee par M. Darboux : D\'eterminer tous les probl\`emes du calcul des variations dans le plan o\`u les solutions sont toutes des droites du plan ; questions qui me semble digne de g\'en\'eralisations f\'econdes et int\'eressantes.]

\end{quote}

\end{document}